\documentclass[10pt]{article}
\usepackage{amssymb}
\usepackage{amsmath}
\usepackage{amsthm}
\usepackage{amsfonts}
\usepackage{amssymb,amsbsy,amsmath,amsfonts,amssymb,amscd}
\def \Prob {{\bf P}}
\newcommand{\R}{\mathbb{R}}

\newcommand{\U}{\mathbb{U}}
\def \H{\mathbb{H}}

\def \expect {{\bf E}}

\def \gamm {{\eta}}
\def \h {{\phi}}

\begin{document}

\title{
Girsanov's transformation for SLE($\kappa, \rho$) processes, intersection exponents and
hiding exponents}

\author{
Wendelin Werner
}

\date {Universit\'e Paris-Sud and IUF}
\maketitle

\begin{abstract}
We relate the formulas giving 
Brownian (and other)
intersection exponents to the absolute 
continuity relations between Bessel process of different
dimensions, via the two-parameter family of Schramm-Loewner Evolution 
processes SLE($\kappa, \rho$)
introduced in \cite {LSWr}. This allows also to
compute the value of some new exponents (``hiding exponents'')
related to SLEs and planar Brownian motions.
\end {abstract}

\section {Introduction}

The value of the intersection exponents between planar Brownian
motions has been derived in the series of papers \cite {LSW1,LSW2,LSW3,LSWa}
using the relationship with the exponents for the
Schramm-Loewner Evolution process with parameter 6 (in short 
SLE$_6$) that can be computed directly.
For instance, if $B^1, \ldots, B^p$ denote $p$ independent planar
Brownian motions started from $p$ fixed different points on the unit
circle, the probability that the $p$ traces $B^1 [0,t],\ldots, B^p[0,t]$
remain disjoint and all stay in the same (fixed) half-plane decays
like $t^{-\tilde\zeta_p /2}$ as $t$ tends to infinity.
The exponent $\tilde \zeta_p$ is called the half-space intersection exponent
between $p$ Brownian motions and it is proved in \cite {LSW1} that
$\tilde \zeta_p = p (2p+1) /3$, as conjectured in \cite {DK}. 

Before SLE allowed to determine the value of these exponents,
it was shown in \cite {LW1} that in order to understand
these Brownian exponents, it is convenient to introduce
``generalized'' Brownian exponents $\tilde \xi_p (a_1, \ldots, a_p)$
that correspond (in the case where all $a_i$'s are integers) to the
decay of the probability of the non-intersection between $p$ unions
of planar Brownian motions in a half-plane containing respectively
$a_1, \ldots , a_p$
paths. For instance $\tilde \xi_p (1, \ldots , 1) = \tilde \zeta_p$.
In particular, one can define the function
$$
U(a)= \lim_{p\to \infty}
\sqrt { \tilde \xi_p (a,a,\ldots, a) / p^2}
$$
and show (this is not mysterious, it is basically a consequence of 
conformal invariance of planar Brownian motion)
that
\begin {equation}
\label {cascade}
 \tilde \xi_p (a_1, \ldots , a_p )
= U^{-1}( U(a_1) + \ldots + U (a_p) ).
\end {equation}
This, combined with the conjectures by Duplantier-Kwon \cite {DK}
for $\tilde \zeta_p$
allowed to predict the value of $U$ and of the generalized
exponents $\tilde \xi$.
Duplantier \cite {Dqg} then observed that this type of equation
can also be viewed as coming from the quantum gravity formalism,
which provided yet another way to predict the exact form of the
function $U$. 

In the paper \cite {LW2}, the relation between the Brownian exponents
and  the exponents for self-avoiding
walks and critical percolation was pointed out. More precisely, a
``universality'' argument was presented that showed that all
conformally invariant models that possess a certain locality
condition must basically have the same exponents i.e. the same
function $U$. This allowed to recover the predictions (see \cite {LW2}
and the references therein) for the
critical exponents for self-avoiding walks or critical percolation
from the above-mentioned prediction for $U$, and conversely to
show that the value of the Brownian exponents would follow from
the computation of the exponents for any other local conformally
invariant object. This is the strategy that was successfully used
in \cite {LSW1}: Show that SLE$_6$ is local, and 
compute its exponents.
The derivation of the SLE$_6$ exponents (in the half-plane) is in
fact a computation related to the (real) Bessel flow. 
This gave the
rigorous proof of the fact that indeed 
$U(x) = \sqrt { x + 1 / 24 } - \sqrt { 1/24}$
as predicted in \cite {LW1}.

In the recent paper \cite {LSWr}, the same basic idea was
developed in a (slightly but essentially) different setting.
There, the family of random sets satisfying the so-called
conformal restriction property is fully described and classified
(the corresponding probability measures are called ``restriction
measures'').
This leads \cite {LSWsaw} to the precise conjecture that SLE$_{8/3}$
is the scaling limit of the half-plane self-avoiding walk. It also proves 
\cite {LSWr} that the boundary of planar Brownian motion, the 
boundary of
the scaling limit of critical percolation cluster interfaces (that Smirnov
\cite {Sm} proved to be indeed corresponding to SLE$_6$) and
the (conjectured) scaling limit of the self-avoiding walk, do not
only have the same exponents but are in fact the same random
object.
The family of restriction measures is parametrized by a positive real
parameter $a$ that can be interpreted as the number of planar
Brownian motions that this restriction measure is equivalent to.
More precisely, when $a$ is  a positive integer, one
can construct the restriction measure with exponent $a$ by considering the union
of $a$ independent Brownian excursions (i.e. in the half-plane,
Brownian motions started
from the origin that are ``conditioned'' to stay forever in the upper half-plane).
This shows that the half-plane intersection exponents $\tilde \xi$
correspond to intersection exponents between restriction measure
samples. 
Note (but this will not be directly relevant here even if it provides one 
additional motivation, since one would wish to also understand the 
relation with the intersection exponents) that the restriction measures are closely related to 
highest-weight representations of some infinite-dimensional Lie algebras (see
\cite {FW}).
 
As shown in \cite {LSWr}, the restriction measures
 (more precisely, their outer boundary) can be
described via variants of SLE$_{8/3}$ called SLE($8/3, \rho$)
(each $\rho$ corresponds to a value of $a$). As we shall briefly recall in the
next section, SLE($8/3, \rho$) is defined as SLE$_{8/3}$ except that the
driving Brownian motion is replaced by a (multiple of) a Bessel
process (actually, it is a little bit more complicated than that).
We shall see in  the present paper 
that with this SLE($\kappa, \rho$) approach, the computation of the intersection exponents
can be interpreted as the standard absolute continuity relations between
Bessel processes of different dimensions (following from
Girsanov's theorem).

This provides the value of various new exponents, some of which describe
probabilities of events
that are associated to planar Brownian motions: For instance,
consider $n+m$ independent Brownian motions in the complex plane that are started 
from $i$, and stopped at their first 
hitting of the line $\{ \Im (z) =R \}$.
What is the probability that
they all stay in the upper half-plane and that none of the
$n$ first Brownian motions contributes to the ``right-hand'' boundary
of the union of the $n+m$ paths restricted to the strip  
$\{\Im (z) \in [1,R] \}$ (i.e. the $n$ paths are hidden from 
$+ \infty$ by the $m$ other paths -- note that this does not imply 
non-intersection between the paths). When $R \to \infty$, the probability that
this happens decays like a negative power of $R$ and the corresponding
exponent is
$$
n+m  + \frac 14 \left( \sqrt {
 24 n + 
(\sqrt { 1 + 24 m}- 3)^2  } - (\sqrt {1 + 24m} - 3 ) 
 \right)
$$
(the $n+m$ part is just corresponding to the fact that the $n+m$ 
paths remain in the upper half-plane). Let us comment that just as for 
the generalized intersection exponents, the values of these
``hiding''  exponents are
rational only for exceptional values of $n,m$. For instance, even the
exponent for $m=n=1$ is the irrational number $(3+\sqrt {7})/2$.
However, for $m=1$ and $n=4$, the exponent is $7$. These ``hiding'' exponents 
do not seem to have appeared before in the theoretical physics literature.

The fact that such exponents can be determined can seem somewhat 
surprising. It is due to the fact that the SLE($8/3, \rho$) approach makes it 
possible to separate the information given by the boundary of the random sets 
(i.e. the law of the exterior boundary of a union of Brownian paths)
from what happens ``in the inside''. An example of 
such facts is the symmetry of the Brownian frontier as described in \cite {LSWr}.

Last but not least, this description not only provides the values of the intersection exponents, but 
it gives directly the law of the paths that 
are conditioned not to intersect. Of course, all this is very closely related to 
the computations of the exponents in \cite {LSW1} as principal eigenvalues of some differential operators, 
and to the corresponding eigenfunction (and the underlying stationary diffusion, for instance
the diffusion conditioned 
to never hit the boundary of the domain), but is simply formulated in terms of these 
SLE($\kappa, \rho$) processes.

The results are not restricted to the $\kappa=8/3$ case.
Hence, one obtains also ``hiding/intersection
exponents'' in the general case. 
In particular, a non-intersection exponent between $p$ SLE$_\kappa$'s (with some 
Brownian loops added in a proper way) turns out to be simply $p(p-1)/\kappa$.  
Recall from \cite {RS}, SLE$_\kappa$ for all $\kappa \in [4,8]$ are supposed to
correspond to the scaling limit of two-dimensional critical statistical physics
models, and that they are conjectured (see \cite {Du2}) to be closely related to the 
SLE($16/ \kappa, \rho$) processes, so that the exponents are relevant in the study of two-dimensional
critical systems.

As explained at the end of the paper, it also gives a new and simple interpretation of the 
``quantum gravity function'' from \cite {KPZ}
predicted by Knizhnik, Polyakov and Zamolodchikov  
that has been used in various forms to predict the values of exponents of two-dimensional critical 
systems by theoretical physicists (see e.g. \cite {Dqg2} and the references therein).
Loosely speaking, the ``quantum
gravity exponent'' (conjectured to correspond to the same system but on a random lattice) is 
just the value $\rho$ that appears in the SLE($\kappa, \rho$).

When writing up this paper, I had basically the choice between on the one hand being sloppy at times, but with 
reasonable heuristic intuition, or giving precise complete statements and proofs that would hide the intuition
behind stochastic calculus considerations and technical setups.
 I deliberately chose the first option, since I believe that the 
gaps left are reasonable.

\medbreak
\noindent
{\bf Acknowledgements.}
I would like to express my deepest gratitude toward Greg Lawler and Oded Schramm. This paper can be viewed
as an addendum to our joint papers, and it could have been
(co)-authored by them as well. I have also benefited from 
stimulating and enlightening discussions with Julien Dub\'edat. 
 
\section {Background}

\subsection {(One-sided) restriction}

We now recall some facts and notation from \cite {LSWr}.
Define the family ${\cal A}$ of closed subsets $A$ of the closed
upper half-plane $\overline \H$ such that
\begin {enumerate}
\item
$\H \setminus A$ is simply connected.
\item
$A$ is bounded and bounded away from the negative half-line.
\end {enumerate}
To each such $A$, associate the conformal mapping $\h_A$ from
$\H \setminus A$ onto $\H$ such that $\h_A(0)=0$ and 
$\h_A(z) \sim z$ when $z \to \infty$.

\medbreak

We say that a closed set $K \subset \overline \H$ is left-filled if
\begin {itemize}
\item
$K$ and $\H \setminus K$ are both simply connected and unbounded
\item
$K \cap \R= \R_-$
\end {itemize}

\medbreak

We say that a random left-filled set satisfies one-sided restriction
if for all $A \in {\cal A}$,
the law of $K$ is identical to that of $\h_A(K)$ given
the event $\{K\cap A = \emptyset\}$.
It is not very difficult (see \cite {LSWr}) to prove that this
implies that for some positive
constant $\alpha$, one has for all $A$,
\begin {equation}
\label {re}
\Prob [ K \cap A = \emptyset ] = \h_A'(0)^\alpha.
\end {equation}
Conversely (see \cite {LSWr}), for each $\alpha>0$, there
exists a unique random right-filled set satisfying this identity.
Its law is called the one-sided restriction measure with exponent
$\alpha$. It can be explicitly constructed using the SLE($8/3,\rho$)
processes that will be described below.

Define now the family ${\cal A}_t$ just as ${\cal A}$ except that
condition 2. is replaced by 
the condition that $A$ is bounded and bounded away from $0$. This
immediately to ``two-sided'' restriction properties
(see \cite {LSWr}) that we shall also use in the
present paper. These are the random sets such that for all 
$A \in {\cal A}_t$, the law of $K$ is is identical to that 
of $\h_A (K)$ conditionally on $\{ K \cap A = \emptyset \}$.
Again, (\ref {re}) has to hold for all $A \in {\cal A}_t$ and some
fixed $\alpha$.
It is proved in \cite {LSWr}
that this can be realized if and only if $\alpha \ge 5/8$. Furthermore,
the only random simple path that satisfies the two-sided
restriction property is SLE$_{8/3}$ for which  the corresponding 
exponent is $\alpha = 5/8$.

\subsection {Absolute continuity relation between Bessel processes}

Suppose that $(X_t, t \ge 0)$ is a Bessel process of dimension $d\ge 1$
started from $x>0$ i.e. (see e.g. \cite {RY} for more details
on the content of this subsection). In other words,
$$
X_t = x + B_t + \int_0^t \frac {d-1}{2X_s}ds
$$
where $(B_t,t \ge 0)$ is a standard one-dimensional Brownian motion started from $0$.
As customary, we will also use the index $\nu$ related to
the dimension $d$ by
$$ d= 2 + 2 \nu .$$
Recall that $X$ is hits the origin if and only if $d <2$.

Suppose for a moment that $d=2$ and $\mu \ge 0$. Then It\^o's formula shows
immediately that
$$
\log X_t = \log x + \int_0^t \frac {dB_s}{X_s} 
$$
is a local martingale.
It is then possible to apply Girsanov's theorem to understand
(for each fixed $t>0$) the behavior of $X$ under the new
probability measure $Q_t$ defined by
$$
dQ_t / dP = (X_t/x)^\mu \exp \left(-\mu^2 \int_0^t\frac { ds }{2 X_s^2} \right)
$$
(it is standard to check that in this particular case, 
the exponential martingale $\exp (\mu \log X_t - \mu^2 \langle \log X 
\rangle_t /2)$ is a plain martingale):
Under this new probability measure,
$$
\tilde B_s := B_s - \mu \int_0^s \frac {ds}{X_s}
$$
for $s \in [0,t]$ is a Brownian motion.
In other words, as
$$
X_s = x +  \tilde B_s + (1/2 + \mu) \int_0^t \frac {ds}{X_s}
$$
it follows that under this new probability measure,
$(X_s, s \in [0,t])$ is a Bessel process of dimension
$d'=2 + 2 \mu$ i.e. of index $\mu$.

Note that the probability measures $Q_t$ are compatible in the sense that
the $Q_T$ probability of any  $\sigma (B_s, s \in [0,t])$ measurable set
is independent of $T>t$. Hence, one can in fact define a
probability measure $Q$ that coincides with $Q_t$ on
$\sigma (B_s, s \in [0,t])$ for all $t$. Under this probability measure 
$Q$, the process $(X_s, s \ge 0)$ is a Bessel process of index $\mu$. 

Conversely, suppose now that $d > 2$ (i.e.
$\nu > 0 $) and define
$$
dQ_t /dP = (X_t/x)^{-\nu} \exp \left( \nu^2 \int_0^t \frac {ds }{ 2 X_s^2} \right).
$$
Then, under the new probability measure $Q$, the process
$X$ is a two-dimensional Bessel process.

Plugging in these two facts together shows that if $X$ is a Bessel process
of dimension $2 +2 \nu \ge 2$,
then under the
probability measure $Q$ that is induced by the probability measures
$Q_t$ defined by
$$
dQ_t/dP
=(X_t/x)^{\mu - \nu} \exp \left( - (\mu^2-\nu^2) \int_0^t \frac {ds }{2 X_s^2} \right)
,$$
the process $X$ is a Bessel process of index $\mu$ (instead of
$\nu$) started from $a$. As we shall see, this relation between $\mu - \nu$ and
the exponent $\mu^2 - \nu^2$ will  basically be the reason for the particular form of the 
critical exponents (i.e. for the function $U$) in our two-dimensional context.

Note that (unless $\mu = \nu$), $Q$ is not absolutely continuous with respect to $P$
(the limiting behavior of $X$ when $t \to \infty$ depends on its
dimension). Similarly, one can let $a$ go to zero, and interpret
heuristically the result as the relation
between Bessel processes of different dimension that are
started from zero. This is not formally true since
$dQ_t/dP$ is not well-defined anymore
($Q_t$ is singular with respect to $P$ because the almost sure 
behavior of $X$ at time $0+$ depends
on its dimension).

\subsection {The SLE($\kappa, \rho$) processes}

We now recall the definition of the SLE ($\kappa, \rho$) processes.
Suppose that $\kappa >0$ and $\rho > -2$.
Let $X$ denote a Bessel process of dimension
$$
d=1 + \frac {2(\rho +2)}{\kappa}
$$
that is started from $X_0= x := a / \sqrt {\kappa} \ge 0$.

Define
$Y = \sqrt {\kappa} X$ and
$$
O_t = \int_0^t \frac {2 ds} { Y_s},
$$
and also
$$
W_t = Y_t + O_t
$$
so that
$$
W_t = a + \sqrt {\kappa}B_t + \int_0^t \frac {\rho ds}{W_s - O_s}
.$$
Then, one defines the Loewner chain $g_t$ with
driving function $W_t$ i.e.
for all $t \ge 0$ and $z$ in the closed upper half-plane $\overline
\H$,
\begin {equation}
\label {loewner}
\partial_t g_t (z)= \frac {2}{g_t (z) - W_t}
\end {equation}
(as long as $g_t (z)$ does not hit $W_t$). For each $t$, $g_t$ is a conformal
map from a domain $H_t$ onto $\H$, where $H_t$ is the set of points $z \in \H$
such that $|g_s (z)- W_s| > 0$ for $s \in [0,t]$. We call this process SLE($\kappa, \rho$).
When $\rho = 0$, this is the (usual) chordal SLE$_\kappa$ process.

We will for the time being assume that $d \ge 2$ (so that $X$ does not 
hit $0$). This means that 
$$
\rho \ge -2 + \frac \kappa 2.
$$
When $\rho =0$, this corresponds to the fact that $\kappa \le 4$.

Suppose that $Y_0 = a  > 0$. Then, for all fixed $t>0$, the
law of $W[0,t]$ is absolutely continuous (even if the density may not be 
bounded) with respect to that
of $\sqrt {\kappa} B_t$, and therefore the law of the Loewner chain
up to time $t$ is absolutely continuous with respect to that
of SLE$_\kappa$. In particular, it is generated by a continuous
curve (see \cite {RS,LSWlesl}).
If $Y_0= 0$, then this does not hold directly, but  one can apply the
same reasoning to the chains
$g_{t_0 + t} \circ g_{t_0}^{-1}$ to deduce that
SLE($\kappa, \rho$) is generated by a continuous curve, that we shall denote by $\gamma$. In other 
words, $g_t$ is the conformal map from the unbounded connected component 
of $\H \setminus \gamma [0,t]$ onto $\H$ that is normalized at infinity 
by $g_t(z) = z + o(1)$.
This curve is
simple if and only if $\kappa \le 4$ (\cite {RS}).

The process
$$
O_t = \int_0^t \frac {2 \ ds}{W_s -O_s}
$$
should be understood as the left-image of the origin
under $g_t$ (when $Y_0 = 0$, then $0$ can correspond to two prime ends
in $g_t^{-1}(\H)$).
The fact the $d \ge 2$ ensures that the left-image of $0$ is never ``swallowed''
by the SLE($\kappa, \rho$) curve, i.e. that the curve never hits the negative half-line.
On the other hand, the SLE($\kappa, \rho$) hits the positive
half-line if $\kappa >4$.

The intuition behind the drift term when $\rho \not= 0$ is the following:
It is a repulsion from the origin if $\rho >0$ or an attractive force
toward the origin if $\rho<0$. The fact that $d \ge 2$ ensures that the
repulsion/attraction is such that the SLE curve never hits the negative half-line:
For instance, when $\kappa = 6$, the repulsion has to be sufficiently strong so that the SLE
does not swallow the origin (i.e. one must have
$\rho \ge 1$). When $\kappa= 8/3$, the process can be attracted toward
zero (i.e. all the values $\rho \ge -2/3$ work) without swallowing it.

Here since $O_0=0$ and $W_0= a$, we say that the SLE($\kappa, \rho$) process is 
started from $(0,a)$. Similarly, for any $o \le w$, one can define an SLE($\kappa, \rho$)
started from $(o,w)$ by translating the SLE started from $(0, w-o)$ by $o$.
 
Note that the SLE($\kappa, \rho$) curve is obtained
deterministically (via the Loewner chain) from the process $Y$ (or $X$).

\subsection {Restriction and SLE($\kappa, \rho$) processes}

In \cite {LSWr}, it is proved that the boundary of the
sample of a one-sided restriction measure of exponent $\eta$ is 
an SLE($8/3, \rho$) process where 
$$ \eta = \bar \eta (8/3,\rho)= \frac {(\rho+2)(3\rho +10)}{32}$$
(here and in the sequel, we will use the bars to indicate that this is a function
and not a parameter).
It is also shown that if one adds (or decorates) an 
SLE$_{\kappa}$ curve with parameter $\kappa \le 8/3$ with a 
Poisson cloud of Brownian loops with intensity $\lambda$, where 
\begin {equation}
\label {lambda}
\lambda = \lambda_\kappa = \bar \lambda ( \kappa)
= \frac { (8 - 3\kappa) ( 6 - \kappa)}{12 \kappa}
\end {equation}
(of course, this depends on the actual normalizing 
factor in the definition of the loop soup), and then ``left-fills'' 
the obtained set, one obtains a sample of a
(one-sided) restriction measure of exponent
$$ \bar \gamm (\kappa, 0) = \frac {6 - \kappa}{2 \kappa}.$$
The same argument is generalized in \cite {Du2}, where it is shown that 
for all $\kappa \le 8/3$, one 
can decorate the SLE($\kappa, \rho$) curve with a Poisson cloud
of Brownian loops of intensity 
$\lambda_\kappa$ and obtain a one-sided restriction measure sample
with exponent
\begin {equation}
\label {bareta} 
\bar \gamm (\kappa, \rho) =
\frac {(\rho +2)(\rho + 6 - \kappa)}{4 \kappa}
\end {equation}
We refer to \cite {LSWr, Du2} for further details.

The Brownian loop decoration procedure can be roughly summarized as follows:
There exists an infinite measure $M$ supported on 
(unrooted) Brownian loops in the half-plane. A realization of the Brownian
loop-soup with intensity $\lambda$
is a Poisson point process with intensity $\lambda M$. A sample of the loop-soup is 
therefore an infinite countable collection of Brownian loops in the upper half-plane. 
One decorates a curve with the loop-soup by adding to the curve all the loops of the
loop-soup that it intersects.
See \cite {LW,LSWr} for more details.
When $\kappa=2$, this is also closely related to the fact that SLE$_2$ is the scaling limit
of loop-erased random walk as proved in \cite {LSWlesl}.

\section {Absolute continuity
between SLE($\kappa, \rho$)'s}
\label {S3}

We are now going to combine the previous considerations.
Consider an SLE($\kappa, \rho$) with $\nu \ge 0$ that is started from
$(0,a)$ as before, where $a>0$.
Recall that
$$
\nu =
\frac {\rho+2}\kappa - \frac 12
$$
i.e.,
$$ 
\rho = \kappa (\nu + \frac 1 2 ) - 2.
$$
Define for $\mu \ge \nu$,
the probability measure $Q$ induced by the measures $Q_t$ as before
i.e. for all $t>0$
\begin {eqnarray*}
 {dQ_t} / {dP}
&=&
 ( X_t / x)^{\mu - \nu}
 \exp \left( \frac { \nu^2 - \mu^2}2 \int_0^t \frac {ds}{ X_s^2 }\right)
 \\
 &=&
 (Y_t/a)^{\mu - \nu} \exp \left( \frac { \kappa (\nu^2 - \mu^2)} 2
 \int_0^t \frac {ds}{Y_s^2} \right)
\end {eqnarray*}
where
$$
X_s = \frac {Y_t} {\sqrt {\kappa}}
= \frac { W_t - O_t} {\sqrt {\kappa} } .
$$
Then, under the probability measure $Q$, the
process $X$ is a Bessel process of dimension $2 + 2 \mu$
(instead of $2+2 \nu$) started from $x = a / \sqrt {\kappa}$.
Hence, under this new probability measure,
the Loewner chain $g_t$ corresponds to that of an
SLE($\kappa, \overline \rho$), where
$$
\bar \rho = \kappa (\mu+ \frac 12)  -2 .
$$
Recall that
$$
\partial_t g_t' (z) = \frac {-2 g_t'(z) }{ (g_t (z) - W_t )^2}
$$
(this formally follows from the differentiation of (\ref {loewner}) with respect to $z$).
Therefore,
$$
\exp \left( \frac {\nu^2 - \mu^2}2 \int_0^t \frac {ds}{ X_s^2} \right)
=
g_t'(0)^{ (\mu^2 - \nu^2)\kappa /4}.
$$
Hence, one can interpret this quantity as the probability that
a sample $K$ of a one-sided restriction measure of exponent
$$
\alpha =
\frac {(\mu^2 - \nu^2) \kappa } 4
$$
avoids $\gamma [0,t]$.
This shows in particular that
\begin {equation}
\Prob [ K \cap \gamma [0,t] = \emptyset ]
=
E[ g_t' (0)^\alpha ] 
= E_{Q_t} [ (x/X_t)^{(\mu - \nu)} ]
= E [ (x/ \tilde X_t)^{(\mu - \nu)}]
\end {equation}
where $\tilde X$ is a Bessel process of dimension $2 + 2 \mu$ started
from $x$.
Let us stress that this is an exact formula and not just an asymptotic expansion.

We now let $a \to 0$ for  fixed $t$. The previous formula shows readily that
$$
\Prob [ K \cap \gamma [0,t] = \emptyset ]
\sim c a^{\mu - \nu},$$
where
$$c = E [ ( \sqrt {\kappa }
\tilde X_t)^{\nu - \mu}]
= (\kappa t)^{(\nu - \mu)/2}
E[ \tilde X_1^{\nu- \mu}]
$$
with $\tilde X_0 = 0$
(the density of $\tilde X_1$ near $0$ behaves like
$y^{1+2\mu}$ so that this expectation is finite).

One can therefore interpret the probability measure $Q$ as follows:
One ``conditions'' $\gamma$ not to intersect a sample of the one-sided 
restriction measure of parameter $\alpha$.
In the same sense as before, this change of measure is singular
if $\mu \not= \nu$, and should be
understood in the $t \to \infty$ (and/or $a \to 0$) limit.
Furthermore, the law of this conditioned SLE($\kappa,\rho$) is that
of SLE($\kappa, \overline \rho$).

To avoid notational confusion and for future 
reference, let us sum up the relation between the exponents, 
$\rho$'s, $\alpha$'s etc.

\begin {itemize}

\item
An SLE$_\kappa$ conditioned  to avoid a sample of a one-sided restriction 
measure of exponent $\alpha$ is an SLE($\kappa, \bar \rho$) where
\begin {equation}
\label {barrho0}
\bar \rho = \bar \rho (\kappa, 0,\alpha) :=
\kappa \sqrt{\frac {4\alpha} \kappa + \left( \frac {2} \kappa
- \frac 12 \right)^2 } + \frac \kappa 2 - 2.
\end {equation}
Conversely, an SLE ($\kappa, \rho$) can be viewed as an SLE$_\kappa$ conditioned 
not to intersect a one-sided restriction sample of exponent
\begin {equation}
\label {baralpha}
\bar \alpha = \bar \alpha (\kappa, \rho) :=
\frac { \rho ( \rho +4 - \kappa)}{4 \kappa}.
\end {equation}

\item
An SLE($\kappa, \rho$) conditioned to avoid a one-sided restriction sample of exponent
$\alpha$ is an SLE($\kappa, \bar \rho$) where 
\begin {equation}
\label {barrho}
\bar \rho = \bar \rho (\kappa, \rho, \alpha)=
\kappa \sqrt{\frac {4\alpha} \kappa + \left( \frac {\rho+2} \kappa
- \frac 12 \right)^2 } + \frac \kappa 2 - 2.
\end {equation}

\item
The exponent associated to the non-intersection event between an SLE($\kappa, \rho$) and 
a one-sided restriction sample of exponent $\alpha$ is
\begin {equation}
\label {barsigma} 
\bar \sigma = \bar\sigma (\kappa, \rho, \alpha)
=
 \sqrt {\frac {4\alpha} \kappa + \left( \frac {\rho+2} \kappa
- \frac 12 \right)^2 }
-
\left( \frac {\rho+2} \kappa
- \frac 12 \right)
. \end {equation}
More precisely, if an SLE($\kappa,\rho$) is started from $a>0$ and runs up to time 1, the probability that 
it does not intersect the one-sided restriction sample of exponent $\alpha$
decays like a constant times $a^\sigma$ when $a \to 0$.
\end {itemize}

Note that 
$$
\bar \rho ( \kappa, \rho, \beta ) 
= \bar \rho ( \kappa, 0, \beta + \overline \alpha ( \kappa, \rho )),
$$
which is not surprising: Conditioning the SLE to avoid a restriction sample of 
exponent $\alpha$ and then to avoid a restriction sample of exponent
$\beta$ is the same as conditioning to avoid a restriction sample of 
exponent $\alpha + \beta$.

Let us briefly insist on the following fact: Here, we focus only on the right-boundaries of the 
domains. For instance, we consider an SLE conditioned not to intersect a one-sided restriction sample, but then 
we focus only on the conditioned SLE and forget about the restriction sample. However, the restriction 
property shows easily how to get the law of the restriction sample conditioned not to intersect that SLE.
More precisely, consider an SLE($\kappa, \rho$) as before. Let $\Gamma_-$ be the connected component of 
$\H \setminus \gamma$ which has the negative half-line on its boundary. Define also a conformal map 
$\Psi_-$ from $\H$ onto $\Gamma_-$ that fixes both the origin and infinity, and let $K$ denote an 
independent sample of the one-sided restriction measure of exponent $\bar \alpha (\kappa, \rho)$.
Then, the joint law of a restriction measure sample conditioned not to intersect an SLE$_\kappa$
(in the sense described above) is just 
that of $(\Psi_- (K), \gamma)$.

\section {Exponents}

This implies a variety 
of results concerning the value of critical exponents.
To illustrate this, we now briefly describe some of them, leaving details 
and further exponents for the interested reader.

\subsection {Hiding exponents between one-sided restriction measures}

Let us first focus on the case $\kappa=8/3$.  In this case,
SLE($8/3, \rho$) is itself the right-boundary of a one-sided
restriction sample of exponent
$\eta = \bar \gamm (8/3, \rho) $.
Suppose that $d \ge 2$ i.e. $\eta \ge 1/3$ and $\rho \ge -2/3$.
Then, the intersection exponent between this right-boundary and another 
one-sided restriction sample of exponent $\beta$ is
\begin {eqnarray*}
\lefteqn{
\sigma = \sigma (\eta \hbox { hides } \beta)
= \bar \sigma ( 8/3, \rho, \beta )}
\\
&=&
 \frac 14
\left( \sqrt {10 +24 \gamm + 24 \beta - 6 \sqrt {1 + 24 \gamm}}
- \sqrt {10 +24\gamm - 6 \sqrt {1+24 \gamm}} \right)
\\
&=&
\frac 14 
\left(
\sqrt { 24 \beta + (\sqrt {1+24 \gamm} - 3)^2 } 
- ( \sqrt {1 + 24 \gamm} - 3 ) \right).
\end {eqnarray*}
This can be interpreted as a hiding exponent between one-sided restriction measures of
exponents $\gamm$ and $\beta$:
Consider two independent samples $K_\gamm$ and $K_\beta$
of one-sided restriction measures with respective exponents $\gamm \ge 1/3$ and
$\beta \ge 0$. Consider the probability that the right-boundary of
$K_\gamm \cup K_\beta$ in the strip $\{\Im (z) \in [1,R] \}$
consists only of points in $K_\gamm$. This probability decays
like $R^{-\sigma}$ as $R \to \infty$. 

In the special case where $\gamm=5/8$, the right-boundary of the 
restriction measure sample is the SLE$_{8/3}$ curve itself.
Hence,  non-intersection between the right-boundary of $K_{5/8}$ and
$K_\beta$ is just non-intersection between $K_\beta$ and the SLE curve,
so that the exponent
$\sigma $ in this case is the same as the non-intersection exponent 
$\hat \xi (5/8, \beta)= \tilde \xi (5/8, \beta) - 5/8 - \beta$ between 
restriction measures.
This gives another way (if one combines the obtained value of 
$\tilde \xi (5/8, \beta)$ with the 
cascade formula (\ref {cascade})) to recover the Brownian
half-plane exponents
that were derived in \cite {LSW1,LSW3} using computations 
involving the SLE$_6$ processes. 

In the very special case where $\eta = \beta = 5/8$, one gets a description of the 
right-boundary of the union of two Brownian excursions in terms of one SLE$_{8/3}$ conditioned 
not to intersect another independent one. We will come back to this is the two-sided case.

When $\gamm>5/8$, the hiding exponent $\sigma$ is smaller than $\hat \xi
(\gamm, \beta)$, which is not surprising since the corresponding events 
are larger.

Note again, that for the values of $\gamm$ such that
$1 + 24 \gamm$ is a perfect square (for instance
$\gamm=1$ and $\gamm=2$ corresponding to one or to the union of two
Brownian excursion), the obtained exponents are
simpler.

Let us insist on the fact that this is 
 valid for all $\beta \ge 0$, and 
$\eta \ge 1/3$ (corresponding to the $d \ge 2$ assumption) but that 
it a priori not holds for $\eta < 1/3$;
we shall see in the next subsection what to do in this case.
 Recall that $\eta =1/3$ corresponds to the scaling limit of conditioned 
percolation cluster boundaries (see \cite {LSWr}).
The exponent $\sigma$ is in this special case equal to $\sqrt {3 \beta/2}$.

There exist various alternative ways to formulate the ``hiding events'' since
restriction measures can be described in terms of Brownian excursions or 
conditioned SLE$_{6}$'s (see \cite {LSWr}).
For example, one way to phrase this in terms of planar Brownian motions is described
at the end of the introduction.
The proof is a consequence of the previous considerations, and of the relation 
between exponents for Brownian excursions and for Brownian motions as developed for instance
in \cite {LSWup}.

Note that the existence of these hiding exponents itself is a 
non-trivial fact (sub-multiplicativity does not simply hold as it does
for non-intersection events).

\subsection {When $\eta < 1/3$}

As we have just mentioned, the previous expression for the hiding exponent $\sigma$ is not 
valid when $\eta < 1/3$. One way to circumvent the difficulty is to first condition the 
boundary of the one-sided restriction sample of exponent $\eta$ not to hit 
the negative half-line (it cannot hide another restriction measure if it hits the negative 
half-line).
Recall that a Bessel process $X$ of dimension $d <2$ started from $x \in [0,1]$ hits $1$ before $0$
with probability $x^{2-d}$ (because $X^{2-d}$ is a local martingale), and that the 
process ``conditioned'' not to intersect $0$ is a Bessel process of dimension $4-d$.

Suppose now that $\rho \in (-2, -2/3)$, so that $d <2$.
The probability that an SLE($8/3,\rho$) started from $(0,a)$ does not intersect 
the negative half-line before its capacity (the Loewner time-parametrization)
 reaches one, decays like a constant times
$a^{2-d}= a^{-1/2 - 3 \rho /4} $ when $a \to 0$.
 Furthermore, the conditioned process is an SLE$(8/3, \rho^*$)
where 
$$ \rho^* = \kappa - 4 - \rho = - \frac 4 3 - \rho$$
(see \cite {Du2} for a similar facts). The corresponding exponents $\eta
= \bar \eta ( 8/3, \rho)$ and $\eta^*= \bar \eta (8/3, \rho^*)$ satisfy
$$ \sqrt { 1 + 24 \eta} + \sqrt {1 + 24 \eta^* } = 6.$$
Straightforward computations then show that for all $\beta >0 $ and $\gamm \in (0, 1/3]$,
\begin {eqnarray*}
 \lefteqn {\sigma ( \gamm \hbox { hides } \beta )}
\\
&=&
\bar \sigma ( \bar \eta ( 8/3, \rho^*)  \hbox { hides } \beta ) 
+ 2 - d \\
&=&
\frac 14 
\left(
\sqrt { 24 \beta + (3-\sqrt {1+24 \gamm})^2 } 
- ( 3 - \sqrt {1 + 24 \gamm}) + ( 6 - 2 \sqrt {1 + 24 \gamm}) \right)
.\end {eqnarray*}
Hence, one can sum up things by saying that the formula
\begin {equation}
{\sigma ( \gamm \hbox { hides } \beta )}
= \frac 14 
\left(
\sqrt { 24 \beta + (3-\sqrt {1+24 \gamm})^2 } 
+ ( 3 - \sqrt {1 + 24 \gamm}) \right)
\end {equation}
 in fact holds 
for all $\gamm >0$ and $\beta >0$.
Let us note that when $ \eta \to 0_+$, one gets a non-trivial limit:
$$
{\sigma (0_+ \hbox { hides } \beta) }
= \frac {1 + \sqrt { 1 + 6 \beta }}2 
,$$
which is somewhat surprising (one might have guessed at first sight that the 
exponent should blow up when $\eta \to 0_+$). Indeed, for each fixed 
large $R$, the probability that $K_{1/N}$ hides $K_{1}$ (with obvious notation) in the 
strip $\{\Im (z) \in [1,R] \}$ is anyway smaller that $1/(N+1)$. This is due to the fact that 
a restriction measure of exponent $(N+1)/N$ can be viewed as the union of $N+1$ independent 
copies of $K_{1/N}$, so that the probability that $K_{1/N}$ hides all $N$ others is no 
larger than $1/(N+1)$. However, when $\eta \to 0_+$ (i.e. $N \to \infty$),
even if the probabilities (for fixed $R$)
go to zero, this does not affect the exponents (only the  
``multiplicative constants'' vanish).

\subsection {Iterations}
The description of conditioned SLE($\kappa, \rho$) as another 
SLE($\kappa, \bar \rho$)
 allows to iterate the procedure, and to obtain exponents describing the 
joint behavior of more than two restriction measures.
For instance, in the simplest case where $\kappa=8/3$, one gets readily the 
exponents describing the non-intersection between $p$ SLE$_{8/3}$'s (these are 
the exponents corresponding \cite {LSWsaw} to the non-intersection of self-avoiding 
walks in a half-plane):

For each positive integer $p$, consider $p$ independent SLE$_{8/3}$'s that are 
conditioned not to intersect (appropriately defined). Define $\eta_p$ the 
exponent of the obtained restriction measure, and define $\rho_p$ 
such that the right-most SLE is an SLE($8/3, \rho_p$). Clearly,
$\eta_p = \bar \eta ( 8/3, \rho_p )$.
Furthermore, for each $p \ge 0$,
$$ \rho_{p+1}  
 =
\bar \rho ( 8/3, 0, \eta_p)
 = 
\bar \rho ( 8/3, 0, \bar \eta ( 8/3, \rho_p ))$$
(where $\rho_1 = 0$).
This shows readily that 
\begin {equation}
\label {etap}
\rho_p = 2 ( p-1) 
\end {equation}
and
\begin {equation}
\eta_p = \frac { p (3p +2)}{8}
.\end {equation}
Hence, the exponent describing the probability that 
$p$ independent chordal SLE$_{8/3}$ (up to time $1$) started at distance $a$ 
of each other are mutually avoiding is 
$$ \eta_p - p \frac 5 8 
= \frac {3p (p-1)} 8 .
$$
This result is not new since (in the notation of \cite {LSW1,LSW3})
$\eta_p = \tilde \xi_p (5/8, \ldots, 5/8)$; these exponents also 
correspond to those conjectured in \cite {DS} for self-avoiding walks
(see \cite {LSWsaw} for the conjectured relation between self-avoiding walks 
and SLE$_{8/3}$).

\subsection {Other $\kappa$'s}

One can easily generalize the iterative procedure for other $\kappa$'s.
Suppose for instance that we consider the conditioned measure for 
$p$ SLE$_\kappa$'s for $\kappa < 8/3$ that are conditioned to mutually avoid each other 
and by the event that no Brownian loop in the Brownian loop-soup with intensity
$\lambda_\kappa$ intersects two different paths. The right-most path is then an 
SLE($\kappa, \rho_p$) for some $\rho_p$ that a priori depends on $\kappa$, but it 
turns out that
$$ 
\rho_p = 2 ( p-1) .
$$
If one adds another independent
Brownian loop-soup with intensity $\lambda_\kappa$ to this right-most path and looks
at the obtained right-most boundary, one obtains a restriction measure with 
exponent
\begin {equation}
\eta_p (\kappa)
=
\bar \eta ( \kappa, \rho_p ) = p  \frac {(2 p + 4 - \kappa)}{2 \kappa} .
\end {equation}
For $\kappa =2$, the exponents correspond to those 
 for loop-erased random walks  
derived by Kenyon \cite {K} and Fomin \cite {Fo} (previously conjectured in \cite {Dle}). This is 
not surprising since loop-erased random walks converge to SLE$_2$ in the 
scaling limit (see \cite {LSWlesl}).

One equivalent way to describe the corresponding event goes as follows:
Run $p$ independent chordal SLE$_{\kappa}$'s 
$S_1, \ldots, S_p$ started from 
nearby points (for instance from the points $a$, $2a$, \ldots, $pa$) up to time one.  
Consider $p$ independent Brownian loop-soups of intensity $\lambda_\kappa$, 
and define for each $j \le p$, the union ${\cal S}_j$ of the loops in the $j$-th 
soup that intersect $S_j$.
Consider now the event that for $j =2$ up to $j=p$,
$$
S_j  \cap ( S_{j-1}  \cup {\cal S}_{j-1} ) = \emptyset
.$$
Then, the probability of this event decays like 
$a^\sigma$ when $a \to 0$, where 
\begin {equation}
\sigma = \eta_p (\kappa) - p \bar \eta ( \kappa , 0 ) 
= 
\frac {p (p-1)}{\kappa}.
\end {equation}
In the special case $\kappa = 2$ that we just mentioned,
 the relation between SLE$_2$ and loop-erased random walks 
\cite {LSWlesl} and Wilson's algorithm \cite {Wi} gives to this event a natural 
interpretation in terms of uniform spanning trees. 

For $\kappa \ge 8/3$, the previous description does not make much 
 sense (the density of 
the loop-soup is negative), and it raises the interesting problem to find a 
simple geometric way to interpret the exponent in terms of a physical model.

\subsection {The ``quantum gravity'' function}

As the formulas show, $\bar \rho ( \kappa, 0, \alpha)$  is in fact the same as the quantum gravity function $U$
(this actually also holds for $\kappa \not= 8/3$), if one compares with the KPZ
equation \cite {KPZ} used e.g. in
Duplantier \cite {Dqg2}. Hence, the SLE($\kappa, \rho$) approach does give another 
interpretation of the ``quantum gravity equations,'' and also permits (using the 
relation with restriction measures) to identify precisely what exponents (i.e. 
what events) are given by this formalism (see Duplantier's \cite {Dqg2}). When $\kappa \not= 8/3$,
this was not so obvious.

On a rigorous level, since the exponents computed via SLE (for instance in the present paper) are 
rigorously derived, while the KPZ relation is not, one may view the SLE derivation of the exponents
as a derivation of the KPZ relation (modulo the assumption that the critical exponents for statistical 
mechanics systems on a random planar graph exist and are universal). 

\subsection {Negative $\alpha$'s}

In fact (but we prefer to emphasize it in this separate paragraph), the absolute 
continuity relation and the derivation of the 
hiding exponents also apply for (some) negative $\alpha$.
In order for the absolute continuity between Bessel processes to hold, the 
condition is that both have a dimension not smaller than $2$. In other words, 
if one starts with an SLE($\kappa, \rho$) such that 
$$\rho \ge -2 + \frac {\kappa}2$$
then, the arguments developed in Section \ref {S3} go through except that there 
is no interpretation of the weighting as a non-intersection probability (the weighting 
is here an unbounded function of the path).
The constraint that the obtained conditioned Bessel process has dimension at least 2 means 
that 
$$ 
\alpha \ge - \frac { (4- \kappa)^2 } {16 \kappa}  
$$
(note that this does not depend on $\rho$). 
Loosely speaking, when $\alpha$ is too negative, then the SLE is not
able to compensate the weighting (so that $Q$ is still a probability measure).
This basically shows that - as one might have expected from the formulas - that 
the hiding exponents make sense on the interval of values of $\alpha$ for which 
it can be extended analytically (as a function of $\alpha$).

In the special case where $\kappa=8/3$, the lower bound on $\alpha$ is $-1/24$.
In the special case $\rho=0$,
the hiding exponent is the intersection exponent 
$\tilde \xi (5/8, \alpha ) - 5/8 - \alpha$. 
We have just argued that $a$ to this power describes indeed 
the asymptotic behavior of the quantity 
$ E [ g_1' (0)^{\alpha} ]$ for an SLE$_{8/3}$ started from $a$ as $a$ vanishes, for all 
values of $\alpha \ge - 1/24$.

If one then applies the cascade ideas, as developed in  \cite {LW1} say, it is then simple to 
see that this for instance enables to deduce that exponents 
 $ \tilde \xi (1, \alpha)-1 - \alpha$ for instance describe the asymptotic behavior
 of $E [g_1' (0)^\alpha ]$ when $g_1$ corresponds this time to the conformal map associated 
 to a Brownian excursion started from $a$, up to time $1$, when $a \to 0$. 
Recall that
$$ 
\tilde \xi (u , \alpha ) = \frac { (\sqrt {24 u +1} + \sqrt {24 \alpha +1 } - 1 )^2 - 1 }{24}. 
$$
In particular,
$$
\tilde \xi (1, -1/24) = 5/8 
\hbox { and } 
\tilde \xi (5/8, -1/24 ) = 1/3.
$$

\section {The two-sided picture}

\subsection {The SLE($\kappa,\rho$) martingales}

Before turning our attention to the two-sided picture, let us point out the 
following  by-product of the description of the SLE($\kappa,\rho$)'s
as an SLE$_\kappa$ conditioned not to intersect a one-sided restriction
sample of exponent
$\alpha (\rho)$. It 
is a simple heuristic explanation to the (useful) martingales
associated to SLE($\kappa,\rho$) derived and
used in \cite {LSWr,Du2}.
Let us first focus on the $\kappa=8/3$ case
studied in \cite {LSWr}.

Let $A \in {\cal A}$. Consider the event that
the SLE($8/3, \rho$) does avoid $A$. Let us now
focus on the conditional probability of this event
given the path up to time $t$.
This is a function of  $W_t$, $O_t$
and of the image of $A$ under $g_t$. Define as in
\cite {LSWr} the conformal map
$h_t$ from $\H \setminus g_t (A)$ onto $\H$
that is normalized by $h_t (z) = z + o (1)$
when $z \to \infty$ (this is just a real shift of $\h_{g_t(A)}$ as
defined in the preliminary section).
If one views the SLE($8/3, \rho$) as an SLE$_{8/3}$
conditioned to avoid a restriction sample $K$, the conditional
probability can be decomposed as follows.
First, the SLE$_{8/3}$ started from $W_t$ has to avoid
$g_t(A)$: This event has probability $h_t'(W_t)^{5/8}$.
Second, the restriction sample has to avoid the
set $g_t(A)$ as well. This occurs with probability
$h_t' (O_t)^\alpha$. Conditionally on these two events,
the image under $h_t$ of the SLE$_{8/3}$ is an SLE$_{8/3}$
in $\H$ started from $h_t (W_t)$ and the image of the
restriction measure sample is a restriction measure sample
in $\H$ started from $h_t (O_t)$. The ``probability'' of
non-intersection between these two sets is going to be affected
by the scaling factor given by the non-intersection exponent
$\nu - 1/4 = \rho/ \kappa $
i.e.
$$
\left( \frac { h_t (W_t) - h_t (O_t)} {W_t - O_t } \right)^{3\rho/8}.
$$
Hence, the quantity 
$$
M_t = 
h_t'(O_t)^\alpha h_t'(W_t)^{5/8} \left( 
\frac { h_t (W_t) - h_t (O_t)} {W_t - O_t } \right)^{3\rho/8} 
$$ 
is a martingale. This is proved analytically in \cite {LSWr}.

The same argument can be used for the local martingales associated to 
SLE($\kappa,\rho$)'s for $\kappa \not= 8/3$ as derived in \cite {Du2}
(with an additional ``loop-soup term'').

\subsection {The two-sided case}
In fact, $M_t$ is a martingale also in the two-sided case. More precisely, 
suppose that $A$ is the symmetric image with respect to the imaginary line 
of a set in ${\cal A}$ (i.e. it is attached to the negative half-line). We 
will suppose in this subsection that $\kappa=8/3$ and $\rho > 0$.  
Then, $M_t$ is still a bounded martingale (this is proved in
\cite {LSWr}), that is well-defined up to the 
(possibly infinite) time $T$ at which the SLE curve hits $A$.
Just as when $A \in {\cal A}$ (see \cite {LSWr}):
\begin {itemize}
\item
If $T$ is finite, then there exists sequence $t_n \to T$
such that $\lim_{n \to \infty} M_{t_n} = 0$.
\item
If $T$ is infinite, then there exists an unbounded sequence $t_n$ 
such that 
$$\lim_{n \to \infty}
h_{t_n}' (W_{t_n}) = \lim_{n \to \infty} \frac 
{ h_{t_n} (W_{t_n}) - h_{t_n} (O_{t_n})} {W_{t_n} - O_{t_n} }
= 1.$$
This is basically due to the fact that $g_t (A)$ becomes smaller and smaller, so that 
$h_t$ becomes closer to the identity.
\end {itemize}
However, the term $h_t' (O_t)$ does not tend to one, because even if $g_t(A)$ becomes smaller, 
$O_t$ gets closer and closer to $g_t (A)$.
But since the SLE path is transient, the term $h_t' (O_t)$ has a (non-trivial) limit when 
$t \to \infty$ (if $T = \infty$) that can be interpreted as follows:

The SLE($8/3, \rho$) is a simple curve $\gamma$ that separates the upper half-plane into two 
connected components $\Gamma_-$ and $\Gamma_+$ (defined in such a way that 
the negative half-line is on the boundary of $\Gamma_-$).
We now focus on $\Gamma_-$. Let $\Phi_-$ denote a anti-conformal map (i.e. $\Phi_- (\bar z)$ is analytic)  
from $\H$ onto $\Gamma_-$ such that $\Phi_- (0) = 0$ and $\Phi_- (\infty) 
= \infty$ (i.e. $\Phi_- (x+iy) = \Psi_- (-x+ iy)$ where $\Psi_-$ is as before).
 In particular, the image of the positive half-line is the negative half-line, and the image of the 
 negative half-line is the curve $\gamma$.
Consider a sample $K$ of a one-sided restriction measure of exponent 
$\alpha$ that is independent of the SLE $\gamma$. Define
$$
{\cal K} = \overline { \Phi_- (K) }.
$$
Note that the set ${\cal K}$ consists of $\gamma$ and of a subset of $\Gamma_-$. In particular,
its ``right-boundary'' is $\gamma$. Since, $K$ is scale-invariant, the actual choice of $\Phi_-$
does not change the law of ${\cal K}$.
Then, the claim is that when $T = \infty$, 
$$
\lim_{t \to \infty}   h_t' (O_t)^\alpha = \Prob [\Phi_- (K)  
\cap A = \emptyset  | \gamma ].
$$
In particular, this implies that almost surely,
$$
\lim_{t \to T} M_t = 1_{T = \infty} \Prob [ {\cal K} \cap A = \emptyset  | \gamma ].
$$
Since the martingale is bounded (by one), the optional stopping theorem shows that
$$
\Prob [ {\cal K} \cap A = \emptyset ] 
= \expect [ M_T ] = \expect [M_0 ] = h_0' (0)^{\eta} = \h_A' (0)^\eta.
$$
But, since $\gamma$ satisfies one-sided restriction with exponent $\eta$,
it follows that in fact 
$$
\Prob [ {\cal K} \cap A = \emptyset ]
= \h_A'(0)^\eta
$$    
for all $A \in {\cal A}_t$.
In other words, ${\cal K}$ is a sample of the two-sided restriction measure of 
exponent $\eta$.

In the special case where $\eta =2$, we see that the restriction measure with exponent 2
corresponds to two SLE$_{8/3}$'s conditioned not to intersect. This is closely related to the 
predictions concerning the scaling limits of self-avoiding polygons \cite {LSWsaw}.
 
\subsection {Two-sided exponents}

{\bf Two-sided hiding.}
This description of the two-sided restriction measures
leads naturally to the following hiding exponents: Consider
 two independent two-sided restriction measure samples
$K_\eta$ and $K_\beta$ of respective exponents
$\eta$ and $\beta$ (where $\eta > 5/8$ and $\beta \ge 5/8$) and focus on the 
event that
$$
K_\beta \cap \{\Im (z)  \in [1,R]\}
\subset K_\eta
.$$
Its probability decays like $R^{-\tau}$ as $R \to \infty$, where 
 \begin {eqnarray*}
\tau &=& \tau( \eta \hbox { hides } \beta ) \\ 
&=& \tilde \xi  ( 5/8, \sigma (\alpha \hbox { hides } \beta ) + \beta + \alpha ) -  \eta - \beta
,\end {eqnarray*}
and the conditional law of $K_\eta$ is the 
restriction measure of exponent $\tau + \eta + \beta$. 
A simple computation yields
\begin {equation}
\tau ( \eta \hbox { hides } \beta )
=
\frac
{\sqrt {24 \beta + (\sqrt {1 +24 \eta} - 6 )^2 } - (\sqrt {1+24 \eta} - 6 )}{2}.
\end {equation}
Note in particular that
\begin {equation}
\label {=5}
\tau (1 \hbox { hides } 1 ) = 3
\hbox { and } 
\tau (2 \hbox { hides } 1 ) = 2 
. \end {equation}
In particular, in both these cases, the exponent of the conditioned restriction measure is $5$.

\medbreak
\noindent
{\bf No cut-points.}
A by-product of these calculations is the exponent that describes the decay of the 
probability that a two-sided restriction measure of exponent $\eta > 5/8$ has no cut-point.
More precisely, when $\eta \in 
(5/8, 35/24)$,
the probability that a sample $K_\eta$ of the two-sided restriction
sample of exponent $\eta$ has no cut-point 
inside the strip $\{ \Im (z)  \in [1,R] \}$ decays like $R^{-\delta(\eta)}$ when 
$R \to \infty$, where
\begin {equation}
\delta (\eta) = 6 - \sqrt { 1 + 24 \eta }
.\end {equation}
Furthermore, the conditional law is 
that of the two-sided restriction measure with exponent 
\begin {equation}
\eta'=  \eta + \delta (\eta) = \eta + 6 - \sqrt {1+ 24 \eta}.
\end {equation}
In other words,
\begin {equation}
\sqrt {1+ 24 \eta} + \sqrt{ 1+ 24 \eta'} = 12.
\end {equation}
The computation goes as follows: The conditioned restriction exponent is 
$$\tilde \xi ( 5/8, \bar \eta (8/3, -4/3 - \rho))$$
where $\rho$ is chosen so that
$\tilde \xi (5/8, \bar \eta (8/3, \rho)) = \eta$.

When $\eta \ge 35/24$, the restriction measure sample has a.s. no cut-point, so that the 
problem is not relevant. 
When $\eta < 5/8$, the two-sided restriction measure does not exist.
When $\eta = 5/8$, then $K_\eta$ is almost surely a simple path, so that 
the probability that it has no cut point in an annulus is $0$.
However, when $\eta  \to 5/8+$, one sees that $\delta$ tends to 2, and that the conditional law
``tends'' to that of a restriction measure of exponent $21/8$, that can therefore be viewed 
as the filling of an SLE$_{8/3}$ conditioned to have no cut-point! Of course, since SLE$_{8/3}$
is a.s. a simple curve, this depends a lot on the limiting procedure used to define this 
conditioned object (here: first replace SLE$_{8/3}$ by a restriction measure of exponent
$5/8 + \epsilon$, then condition it to have no cut point (in larger and larger annuli), and 
finally let $\epsilon$ tend to zero).

It is worthwhile stressing the special case where $\eta=1$.   
The exponent $\delta$ is equal to  $1$ and it 
is related to B\'alint Vir\'ag's Brownian beads \cite {V}.
It gives a description 
of the restriction measure of parameter 
$\eta'=2$ as the filling of one single path. More precisely:
``The filling of a Brownian excursion conditioned to have no cut point
has the same law as the filling of the union of two Brownian excursions.''
It raises the question whether 
this conditioned Brownian excursion
has something to do with the path that is obtained by considering the appropriate
SLE$_\kappa$ to which one chronologically attaches Brownian loops as in \cite {LSWr}
in order to construct a restriction measure sample of exponent $2$.

Note also that the two-sided measure obtained if one conditions $K_\eta$ to hide $K_\beta$, is the
same as the one obtained if one conditions $K_{\eta'}$ to hide $K_\beta$. This is not surprising: One first conditions 
$K_\eta$ to have no cut point, and then weights it by the ``space'' it leaves in its inside.
 
\medbreak
\noindent
{\bf Mixed two-sided hiding.}
One can also define exponents associated to ``mixed'' 
two-sided hiding: Consider the exponent $\hat \tau (\eta, \beta)$ that is associated to the 
fact that the left-boundary of $K_\eta \cup K_\beta$ consists only of points in $K_\eta$ 
while the right-boundary consists of points in $K_\beta$. This time
$$ 
\hat \tau 
=
\tilde \xi (5/8, \sigma (\beta \hbox { hides } \alpha) + \alpha + \beta ) - \eta - \beta
,$$
where as before
$\tilde \xi (5/8, \alpha) = \eta$.
This leads to 
\begin {equation}
\hat \tau ( \eta, \beta )
= 
\frac { 9 - B - E + 2 \sqrt {(B-3)^2 + (E-3)^2 -1 }} {4}
\end {equation}
where $B= \sqrt {1 + 24 \beta}$ and $E = \sqrt {1+24 \eta}$.
For instance $\hat \tau (1,1) = (2 \sqrt 7 - 1)/4$.

\medbreak
\noindent
{\bf Radial hiding.}
All two-sided hiding exponents yield
readily the corresponding exponent in the radial setting, using the mapping 
described for example in \cite {LW2} and the disconnection exponents computed in 
\cite {LSW2,LSW3} (see also \cite {LSWr2}). 

For instance, consider $n+p$ independent Brownian motions started from the origin 
and stopped when they hit the unit circle. Consider the event that the union of these 
$n+p$ paths do not disconnect the circle of radius $r$ from $1$, and that the boundary of the 
connected component of $\U \setminus (B^1 \cup \ldots B^{n+p})$ that contains $1$ 
contains no point of $B^{n+1} \cup \ldots \cup B^{n+p}$.
 Then, the probability of this event decays like $r^{\rho}$ when $r \to 0$, where 
\begin {equation}
\rho = \rho ( n \hbox { hides } p )
=
\frac {
\left( \sqrt {24 p + ( \sqrt {1+ 24 n } - 6 )^2 } + 5 \right)^2 - 4 } {48}
\end {equation}
Again, note that $\rho ( \eta \hbox {hides} \beta )
= \rho (\eta' \hbox { hides }\beta )$.
Note also that when $n= 2$ or $n=1$, the hiding exponent is just
$$
\rho ( 2 \hbox { hides } \beta ) = \rho (1 \hbox { hides } \beta )
 =  \xi (2, \beta)
$$ 
(in the notation of \cite {LSW2}), which is 
not surprising because of the inside/outside symmetry of the Brownian frontier 
pointed out in \cite {LSWr}.
The inside/outside symmetry of the Brownian frontier also shows that a single 
Brownian motion started from the origin, ``conditioned not to disconnect the origin from 
infinity and to have no cut point'' also separates the plane into the ``inside'' $I$ and the ``outside''
$O$ in such a way that $(I,O)$ and $(O,I)$ have the same law.

When the half-plane exponent is $5$, then the radial exponent is $2$. For instance, 
$\rho (1 \hbox { hides } 1) = 2$, so that the corresponding existence problem is ``critical'': 
Are there points $B_T$ on the outer boundary of a
planar  Brownian path $(B_t, t \in [0,1])$ 
such that (locally) the outer boundary consists only of the future after $B_T$ (or only of the past
before $B_T$)?

------------------

Laboratoire de Math\'ematiques

Universit\'e Paris-Sud

91405 Orsay cedex, France

wendelin.werner@math.u-psud.fr

\end {document}